# Multistability in Monotone I/O Systems, Preliminary Report

David Angeli
Dip. Sistemi e Informatica
University of Florence, 50139 Firenze, Italy
angeli@dsi.unifi.it

Eduardo D. Sontag
Dept. of Mathematics
Rutgers University, NJ, USA
sontag@hilbert.rutgers.edu

## 1 Introduction

We extend the setup in [1] to deal with the case in which more than one steady state may exist in feedback configurations. This provides a foundation for the analysis of multi-stability and hysteresis behavior in high dimensional feedback systems.

We do not repeat here most of the basic definitions and notations from [1], save for the reminder that, with respect to appropriate positivity cones, we say that a system is monotone if:

$$\xi_1 \succeq \xi_2 \quad \& \quad u_1 \succeq u_2 \Rightarrow x(t, \xi_1, u_1) \succeq x(t, \xi_2, u_2) \quad \forall\, t \geq 0.$$

If $\operatorname{int}(K) \neq \emptyset$, this is equivalent to asking:

$$\xi_1 \gg \xi_2 \quad \& \quad u_1 \succeq u_2 \Rightarrow x(t, \xi_1, u_1) \gg x(t, \xi_2, u_2) \quad \forall\, t \geq 0$$

(a set which is the closure of its interior is invariant iff its interior is invariant). We also recall the following definition: a system is strongly monotone if:

$$\xi_1 \succ \xi_2 \quad \& \quad u_1 \succeq u_2 \Rightarrow x(t, \xi_1, u_1) \gg x(t, \xi_2, u_2) \quad \forall\, t > 0.$$

It is often convenient to assume more about the steady-state convergence properties of a monotone system. The following notion, first introduced in [1], will be useful in order to state our main result.

**Definition 1.1** We say that a system admits a *non-degenerate* I/S static characteristic $k_x(\cdot): \mathcal{U} \to X$ if for all constant inputs $u \in \mathcal{U}$ there exists a unique globally asymptotically stable equilibrium $k_x(u)$ and $\det(D_x f(k_x(u), u)) \neq 0$. □

Notice that, for technical reasons, the property has been strengthened with respect to the definition in [1] by assuming non-degeneracy of the equilibria.

## 2 Sufficient Conditions for Strong Monotonicity

Detecting if a system is monotone with respect to the partial order induced by some positivity cone $K$, without actually having to compute explicit trajectories of the system itself, is of course a very important task in order to apply our results in any specific situation. Necessary and sufficient differential characterizations of monotonicity are discussed in [1], where extensions to systems with inputs and outputs are presented of some well-known criteria previously only formulated for autonomous differential equations (see [13]). For the sake of completeness we recall the differential characterization proved in [1]:

**Theorem 1** *A finite-dimensional nonlinear systems of differential equations $\dot{x} = f(x,u)$ with state-space $X$ and input-space $\mathcal{U}$ is monotone with respect to the positivity cone $K$ and $K^{\mathcal{U}}$ for inputs if and only if:*

$$x_1 \succeq x_2 \text{ and } u_1 \succeq u_2 \Rightarrow f(x_1, u_1) - f(x_2, u_2) \in TC(K_{x_1-x_2}) \tag{1}$$

∎

For the special case of positivity orthants, criteria are better formulated in terms of the *incidence graph* of the system. Along similar lines as in [10], given a system:

$$\begin{aligned} \dot{x} &= f(x,u) \\ y &= h(x) \end{aligned} \tag{2}$$

with $x \in X \subset \mathbb{R}^n$, $u \in \mathcal{U} \subset \mathbb{R}^m$ and $h(\cdot) : X \to \mathcal{Y} \subset \mathbb{R}^p$, we associate to it a signed digraph, with vertices $x_1, x_2 \ldots x_n$, $u_1, u_2, \ldots u_m$, $y_1, y_2 \ldots y_n$ and edges constructed according to the following set of rules:

**Edges between $x$ vertices:**

The graph is defined only for systems so that for any couple $1 \leq i, j \leq n$ of integers with $i \neq j$ one of the following rules apply:

1. If $f^i(x,u)$ is strictly increasing with respect to $x_j$ for all $x, u \in X \times \mathcal{U}$ then we draw a positive edge $e_{ij}$ directed from vertex $x_j$ to $x_i$.

2. If $f^i(x,u)$ is strictly decreasing as a function of $x_j$ for all $x, u \in X \times \mathcal{U}$ then we draw a negative edge $e_{ij}$ directed from vertex $v_j$ to $v_i$.

3. Otherwise, $\frac{\partial f_i}{\partial x_j} = 0$ for all $x, u$ and no edge from $x_j$ to $x_i$ is drawn.

**Edges between $u$ and $x$ vertices:**

The graph is defined only for systems so that for any couple of integers $i, j$ with $1 \leq i \leq n$ and $1 \leq j \leq m$ one of the following rules apply:

1. If $f^i(x,u)$ is strictly increasing as a function of $u_j$ for all $x, u \in X \times \mathcal{U}$ then we draw a positive edge $\tilde{e}_{ij}$ directed from vertex $u_j$ to $x_i$.

2. If $f^i(x,u)$ is strictly decreasing as a function of $u_j$ for all $x, u \in X \times \mathcal{U}$ then we draw a negative edge $\tilde{e}_{ij}$ directed from vertex $u_j$ to $x_i$.

3. Otherwise $\frac{\partial f_i}{\partial u_j} = 0$ for all $x, u$ and no edge from $u_j$ to $x_i$ is drawn.

**Edges between $x$ and $y$ vertices:**

The graph is defined only for systems so that for any couple of integers $i, j$ with $1 \leq i \leq p$ and $1 \leq j \leq n$ one of the following rules apply:

1. If $h^i(x)$ is strictly increasing as a function of $x_j$ for all $x \in X$ then we draw a positive edge $\tilde{e}_{ij}$ directed from vertex $x_j$ to $y_i$.

2. If $h^i(x)$ is strictly decreasing as a function of $x_j$ for all $x \in X$ then we draw a negative edge $\tilde{e}_{ij}$ directed from vertex $x_j$ to $y_i$.

3. Otherwise, $\frac{\partial h_i}{\partial x_j} = 0$ for all $x \in X$ and no edge from $x_j$ to $y_i$ is drawn.

Under this convention, a directed path $\mathcal{P}$ is a finite sequence of vertices, $v_{n_0}, v_{n_1} \ldots v_{n_L}$, such that each vertex appears at most once in the sequence and $e_{ij}$ is an edge whenever $v_i, v_j$ appear consecutively in the path. The integer $L$, is called the length of the path and it is denoted by $L(\mathcal{P})$. By $\mathcal{P}_i$, we denote the $v_{n_i}$, the $i+1$-th, vertex of the path $\mathcal{P}$. A *cycle*, not necessarily directed, is a finite sequence of vertices $v_{n_0}, v_{n_1} \ldots v_{n_L}$ such that $v_{n_0} = v_{n_L}$ and the constraint that either $e_{ij}$ or $e_{ji}$ is an edge whenever $v_i$ and $v_j$ appear consecutively in the cycle. The sign of a cycle is the product of the signs of the edges comprising it (accordingly the sign of a path is the product of the signs of its edges). One of the main results in [10] is that an autonomous system (no inputs) is monotone with respect to some orthant if and only if its associated graph does not contain any *negative (simple) cycles*. An analogous result holds for controlled systems:

**Proposition 2.1** A system (2) which admits an incidence graph according to the above set of rules is monotone with respect to some orthants $K$, $K^{\mathcal{U}}$ and $K^{\mathcal{Y}}$ if and only its graph does not contain any *negative (simple) cycles*. □

**Remark 2.2** We remark that in this set-up we deliberately restricted the class of systems for which the incidence graph is defined. In [10] in fact the milder requirement that $\frac{\partial f_i}{\partial x_j} \geq 0$ for all $x$ together with $\frac{\partial f_i}{\partial x_j} > 0$ for some $x$ is asked for in order to draw an edge between vertices $x_i$, $x_j$. This more general notion of incidence graph is however much more cumbersome to deal with if we want to give conditions for strong monotonicity of a system. □

This definition of incidence graph also provides the right set-up for easy geometrical characterizations of the property defined below (see [6] for the linear systems case).

**Definition 2.3** A MIMO system is (weakly) *excitable* if for any initial condition $\xi$ and any couple of inputs $v, u$ with $v(t) \succ u(t)$ ($v(t) \gg u(t)$) for almost all $t > 0$, the following holds:

$$x(t, \xi, v) \gg x(t, \xi, u) \qquad \forall t > 0. \tag{3}$$

□

It turns out that this property can be easily characterized in terms of the incidence graph (see [12]):

**Theorem 2** *A monotone system which admits an incidence graph is (weakly) excitable provided that each $x_i$ is reachable through a directed path from any (some) $u_j$.*

It is worth pointing out that for the case of positive linear systems the above conditions are necessary and sufficient (see [6]).

The proof of this result is based on the following simple Lemma:

**Lemma 2.4** *A scalar monotone system $\dot{x} = f(x, u)$ admitting an incidence graph is (weakly) excitable if and only if, for all (for some) $j$, $f(x, u)$ is strictly monotone as a function of $u_j$.*

*Proof.* Consider, without loss of generality, the case of a system monotone with respect to the standard positive orthants. Let $v(t) \succ u(t)$ be arbitrary input signals. In particular, there exists an integer-valued function of time $j(t)$ so that $v_{j(t)}(t) > u_{j(t)}(t)$ for almost all $t > 0$. By monotonicity we know that $x(t, \xi, v) \geq x(t, \xi, u)$ for all $t \geq 0$. We need to show that the inequality is strict for all $t > 0$. Assume by contradiction[1] that $x(t, \xi, v) = x(t, \xi, u) := x(t)$ for all $t \in (0, \varepsilon)$ for some $\varepsilon > 0$. Then taking derivatives with respect to time we have, $f(x(t), v(t)) = f(x(t), u(t))$ for almost all $t \in (0, \varepsilon)$ but this contradicts strict monotonicity of $f(x, u)$ with respect to all of the $u_j$s.

Conversely, let the system admit an incidence graph. Hence, $f(x, u)$ is either strictly monotone with respect to $u_j$ or $\frac{\partial f}{\partial u_j}(x, u) \equiv 0$ for all $(x, u)$. Assume by contradiction that there exists $\frac{\partial f}{\partial u_j}(x, u) \equiv 0$ for some $j$. Then, given any initial condition $\xi$, and any couple of feasible input values $u, v$ so that $v - u = |v - u| e_j$, we have: $x(t, \xi, v) = x(t, \xi, u)$, for all $t \geq 0$, as the $j$-th component of the input does not influence the solution of the differential equation. But this contradicts excitability.

Similar arguments apply to the case of weak excitability. We give the sketch of the proof below: Consider arbitrary input signals $v$ and $u$ so that $v(t) \gg u(t)$ for almost all $t \geq 0$. By monotonicity $x(t, \xi, v) \geq x(t, \xi, u)$ for all $t \geq 0$. Hence, we only need to show that the inequality is strict for all $t > 0$. Assume by contradiction $x(t, \xi, v) = x(t, \xi, u) := x(t)$ for all $t \in (0, \varepsilon)$ for some $\varepsilon > 0$. Then taking derivatives with respect to time we have, $f(x(t), v(t)) = f(x(t), u(t))$ for all $t \in (0, \varepsilon)$ but this contradicts strict monotonicity of $f(x, u)$ with respect to at least some of the $u_j$s.

The converse implication is trivial as $\frac{\partial f}{\partial u_j} \equiv 0$ for all $j$ implies that solutions do not depend upon input signals, and this clearly violates weak excitability. ∎

Theorem 2 can be proved by induction by applying repeatedly Lemma 2.4 (see appendix for a detailed proof).

The dual of excitability is also useful in the following discussion:

**Definition 2.5** A system is (weakly) *transparent* if for each pair of solution $x(t, \xi_1, u)$, $x(t, \xi_2, u)$ with $\xi_1 \succ \xi_2$ we have $h(x(t, \xi_1, u)) \gg (\succ) h(x(t, \xi_2, u))$ for all $t > 0$. □

**Theorem 3** *A monotone system which admits an incidence graph is (weakly) transparent provided that directed paths exist from any $x_j$ to any (at least one) output vertex $y_i$.*

The proof of Theorem 3 is given in an Appendix. In the present section we discuss a sufficient condition for strong monotonicity of MIMO systems in unitary feedback.

**Theorem 4** *Consider the unitary feedback interconnection of a MIMO system:*

$$\begin{align} \dot{x} &= f(x, u) \\ y &= h(x) \end{align} \quad (4)$$

*viz. the differential equation resulting from (4) when we let $u = y$. The induced flow is* strongly *monotone provided that (4) be monotone, excitable and transparent with either excitability or transparency possibly holding in a weak sense.*

*Proof.* By Theorem 1, we know that

$$x_1 \succeq x_2 \ \& \ u_1 \geq u_2 \ \Rightarrow \ f(x_1, u_1) - f(x_2, u_2) \in TC_{x_1 - x_2}(K) \quad (5)$$

where $K$ is the positivity cone relative to the order $\succeq$. Let us first show monotonicity of the feedback loop system. Recall that $h$ is a monotone map, viz.:

$$x_1 \succeq x_2 \ \Rightarrow \ h(x_1) \geq h(x_2) \quad (6)$$

---
[1] We exploit the fact that the interior of the positivity cone is invariant

Therefore, if we combine (5) with (6) and we let $u_1 = h(x_1)$ and $u_2 = h(x_2)$ we obtain

$$x_1 \succeq x_2 \Rightarrow f(x_1, h(x_1)) - f(x_2, h(x_2)) \in TC_{x_1-x_2}(K) \tag{7}$$

which, by Theorem 1 in [1] is equivalent to monotonicity of the closed-loop system:

$$\dot{z} = f(z, h(z)). \tag{8}$$

In particular then, if we denote by $z(t, \xi)$ the solutions of (8) we have as a consequence of monotonicity:

$$\xi_1 \succeq \xi_2 \Rightarrow h(z(t, \xi_1)) \geq h(z(t, \xi_2)) \qquad \forall\, t \geq 0. \tag{9}$$

Exploiting the fact that $z(t, \xi) = x(t, \xi, h(z(t, \xi)))$ and (weak) strong transparency of (4) we obtain:

$$\begin{aligned}\xi_1 \succ \xi_2 \Rightarrow\ & h(z(t, \xi_1)) = h(x(t, \xi_1, h(z(t, \xi_1)))) \gg (\succ) h(x(t, \xi_2, h(z(t, \xi_1)))) \\ \succeq\ & h(x(t, \xi_2, h(z(t, \xi_2)))) = h(z(t, \xi_2)) \qquad \forall\, t > 0.\end{aligned} \tag{10}$$

Finally, by weak (strong) excitability and (10)

$$\begin{aligned}\xi_1 \succ \xi_2 \Rightarrow\ & h(z(t, \xi_1)) \gg (\succ) h(z(t, \xi_2)) \\ \Rightarrow\ & z(t, \xi_1) = x(t, \xi_1, h(z(t, \xi_1))) \gg x(t, \xi_2, h(z(t, \xi_2))) = z(t, \xi_2) \qquad \forall\, t > 0\end{aligned} \tag{11}$$

as desired. ∎

## 3  Monotone Linear Systems

We recall next some basic facts about linear monotone systems which will be of interest in the discussion of the main result.

**Theorem 5** *Let us consider the following finite dimensional MIMO linear system:*

$$\dot{x} = Ax + Bu, \qquad y = Cx. \tag{12}$$

*with $x \in (\mathbb{R}^n, \succeq_x)$, $u \in (\mathbb{R}^m, \succeq_u)$, $y \in (\mathbb{R}^p, \succeq_y)$ (viz. we assume the state, input and output space equipped with some partial orders induced by the positivity cones $K^{\mathcal{X}}$, $K^{\mathcal{U}}$ and $K^{\mathcal{Y}}$ respectively).*

*System (12) is a monotone control systems with respect to the partial orders specified above if and only if:*

1. *$K^{\mathcal{X}}$ is positively invariant for the autonomous system $\dot{x} = Ax$;*

2. *$BK^{\mathcal{U}} \subseteq K^{\mathcal{X}}$;*

3. *$CK^{\mathcal{X}} \subseteq K^{\mathcal{Y}}$.*

*Proof.* By the characterization of monotonicity in Theorem 1, a system is monotone if and only if:

$$x_1 \succeq_x x_2\ \&\ u_1 \succeq_u u_2 \quad \Rightarrow \quad A(x_1 - x_2) + B(u_1 - u_2) \in TC_{x_1-x_2}(K^{\mathcal{X}}), \tag{13}$$

and the output map is monotone, viz.:

$$x_1 \succeq x_2 \Rightarrow Cx_1 \succeq Cx_2. \tag{14}$$

In terms of positivity cones and denoting $\tilde{x} := x_1 - x_2$ and $\tilde{u} = u_1 - u_2$, conditions (13) and (14) are equivalent to:
$$\tilde{x} \in K^X \ \& \ \tilde{u} \in K^{\mathcal{U}} \ \Rightarrow \ A\tilde{x} + B\tilde{u} \in TC_{\tilde{x}}(K^X) \tag{15}$$

and:
$$\tilde{x} \in K^X \ \Rightarrow \ C\tilde{x} \in K^{\mathcal{Y}}. \tag{16}$$

Condition (13) is clearly equivalent to assumption 3). Condition (15) can be further decomposed by first taking arbitrary $\tilde{x}$ and fixing $\tilde{u} = 0$ and then $\tilde{x} = 0$ and arbitrary $\tilde{u}$. Condition (15) therefore implies (and is in fact equivalent to as we shall see later):
$$\tilde{x} \in K^X \ \Rightarrow \ A\tilde{x} \in TC_{\tilde{x}}(K^X) \tag{17}$$

and:
$$\tilde{u} \in K^{\mathcal{U}} \ \Rightarrow \ B\tilde{u} \in TC_0(K^X) = K^X. \tag{18}$$

The converse implication just follows by recalling that tangent cones of a convex set are closed under sums (since they are convex cones) and the following inclusion holds: $K^X \subseteq TC_{\tilde{x}}(K^X)$ for any $\tilde{x} \in K^X$. Condition (18) is clearly assumption 2). Whereas condition (17) is the well-known characterization of positive invariance of $K^X$ under the flow $\dot{x} = Ax$. ∎

**Corollary 3.1** *The impulse response of a finite-dimensional, monotone, linear system (with respect to positive impulses) is a positive signal in output space:*
$$Ce^{At}BK^{\mathcal{U}} \ \subseteq \ K^{\mathcal{Y}}$$

□

The following fact is a straightforward consequence of the Perron-Frobenius (Krein-Rutman) Theorem (see [2] pp. 6-8):

**Theorem 6** *Assume that the linear system $\dot{x} = Ax$ admits a positively invariant convex (and proper) cone $K$. Then, there exists a dominant real eigenvalue $\lambda$ (viz. an eigenvalue so that $\text{Re}[\lambda_i] \leq \lambda$ for all $i \in 1, 2, \ldots n$), and the corresponding positive eigenvector $v_\lambda$ (unique up to a positive multiple) satisfies $v_\lambda \in K$.*

*Proof.* Consider the exponential map $\xi \to e^{At}\xi$. By positive invariance of $K$, for each $t > 0$ the exponential is a linear map from $K$ to $K$. Moreover, for $t$ sufficiently small it is one-to-one on the spectrum of $A$. Thus, by Lemma A.3.3 in [14], the geometric multiplicity of $e^{\lambda_i t}$ as an eigenvalue of the exponential map is the same as $\lambda_i$ as an eigenvalue of $A$. Therefore, we can study the spectrum of $A$ by looking at the spectrum of its exponential map for $t$ sufficiently small. By the Perron-Frobenius Theorem, there exists a real positive eigenvalue $\mu$, with eigenvector $v \in K$, which is dominant in the sense that $\mu = \rho(e^{At})$ (eigenvalue of maximum modulus). Therefore, we conclude that $\lambda := \log(\mu)/t$ is an eigenvalue for $A$, relative to the same eigenvector $v \in K$, and $\text{Re}(\lambda) \geq \text{Re}(\lambda_i)$ for all $\lambda_i \in \text{Spec}(A)$. ∎

**Remark 3.2** It is worth pointing out that for asymptotically stable SISO monotone systems, the condition $h(t) \geq 0$, implies that the $\mathcal{L}^\infty \to \mathcal{L}^\infty$ induced gain, equals the steady state gain. Recall that the steady-state gain of a linear system is just the slope of its I/O static characteristic. The $L_\infty \to L_\infty$ induced gain is instead defined as:
$$\gamma_\infty := \sup_{u \neq 0} \frac{\|y\|_\infty}{\|u\|_\infty}$$

where $y(t) = y(t, 0, u)$. It is well known (see [5]) that $\gamma_\infty$ equals the $L_1$ norm of the impulse response. Thus,
$$\gamma_\infty = \int_0^{+\infty} |h(t)|\, dt = \int_0^{+\infty} h(t)\, dt = -CA^{-1}B = k'_y(u).$$

□

The next technical lemma will be useful in order to study nonlinear monotone systems by linearizing the flow around an equilibrium position:

**Lemma 3.3** Let $f : X \times \mathcal{U} \to \mathbb{R}^n$ be a $\mathcal{C}^1$ vector-field. Let $f(\bar{x}, \bar{u}) = 0$ for some $\bar{x} \in X$ and $\bar{u} \in \mathcal{U}$. If the flow induced by $f$ is monotone with respect to some positivity cone $K$, the same holds true for the linearization at $(\bar{x}, \bar{u})$:
$$\begin{aligned}
\dot{z} &= \left.\frac{\partial f}{\partial x}\right|_{x=\bar{x}, u=\bar{u}} z + \left.\frac{\partial f}{\partial u}\right|_{x=\bar{x}, u=\bar{u}} v \\
w &= \left.\frac{\partial h}{\partial x}\right|_{x=\bar{x}} z
\end{aligned} \tag{19}$$

*Proof.* By one the results in [1], a system is monotone with respect to the positivity cones $K$ (for states) and $K^\mathcal{U}$ (for inputs) if and only if:
$$x_1 \succeq x_2,\ u_1 \succeq u_2 \Rightarrow f(x_1, u_1) - f(x_2, u_2) \in TC_{x_1 - x_2}(K). \tag{20}$$

Let $z \in K$, $v \in K^\mathcal{U}$ be arbitrary and, for any $\varepsilon > 0$, $x_\varepsilon := \varepsilon z + \bar{x}$, $u_\varepsilon = \varepsilon v + \bar{u}$. By (20) applied with $x_1 = x_\varepsilon$ and $x_2 = \bar{x}$,
$$f(x_\varepsilon, u_\varepsilon)/\varepsilon \in TC_{\varepsilon z}(K) = TC_z(K). \tag{21}$$

By letting $\varepsilon$ tend to 0 and exploiting closedness of the tangent cone we have:
$$z \succeq 0,\ v \succeq 0 \Rightarrow \left.\frac{\partial f}{\partial x}\right|_{x=\bar{x}, u=\bar{u}} z + \left.\frac{\partial f}{\partial u}\right|_{x=\bar{x}, u=\bar{u}} v \in TC_z(K). \tag{22}$$

Let, for simplicity $A = \left.\frac{\partial f}{\partial x}\right|_{x=\bar{x}, u=\bar{u}}$ and $B = \left.\frac{\partial f}{\partial u}\right|_{x=\bar{x}, u=\bar{u}}$. By linearity, there follows easily:
$$z_1 \succeq z_2,\ v_1 \succeq v_2 \Rightarrow (Az_1 + Bv_1) - (Az_2 + Bv_2) \in TC_{z_1 - z_2}(K). \tag{23}$$

This concludes the proof of the claim, by exploiting once more the characterization of monotonicity in [1]. ■

**Lemma 3.4** Consider a monotone system with a non-degenerate I/S static characteristic $k_x(\cdot)$. For each $u \in \mathcal{U}$ the corresponding equilibrium $k_x(u)$ is hyperbolic.

*Proof.* By Lemma 3.3 the linearized system at the equilibrium is monotone. Therefore it admits a real dominant eigenvalue $\lambda$. By asymptotic stability of the nonlinear system and non-degeneracy, $\lambda < 0$. Thus for all $\lambda_i \in \text{sp}[D_x f(k_x(u), u)]$ we have $\text{Re}[\lambda_i] \leq \lambda < 0$ which completes the proof of our claim. ■

We remark that for the special case of $K$, $K^\mathcal{U}$ being positive orthants the result was already proved in Section 8, [1].

# 4 Main Result

Our main result will provide a global analysis tool for positive unitary feedback interconnection of monotone systems. The fixed points of the I/O characteristic will play a central role in the statement of the result. In particular, we say that $k_y(\cdot) : \mathcal{U} \to \mathcal{Y}$ has non-degenerate fixed points if for all $u \in \mathcal{U}$ with $k_y(u) = u$ we have that $k'_y(u)$ exists and $k'_y(u) \neq 1$.

**Theorem 7** *Consider a strongly monotone, SISO, dynamical system, endowed with a non-degenerate I/S and I/O static characteristic:*
$$\begin{aligned} \dot{x} &= f(x,u) \\ y &= h(x), \end{aligned} \qquad (24)$$
*where $f$ is of class $\mathcal{C}^1$. Consider the unitary positive feedback interconnection $u = y$. Then the equilibria are in 1-1 correspondence with the fixed points of the I/O characteristic. Moreover, if $k_y$ has non-degenerate fixed points, and all trajectories are bounded, then for almost all initial conditions, solutions converge to the set of equilibria of (24) corresponding to inputs for which $k'_y(u) < 1$.*

*Proof.* Let $k_x : \mathcal{U} \to X$ denote the I/S static characteristic and $\bar{u}$ any solution of $u = h(k_x(u))$. Clearly, $f(k_x(\bar{u}), h(k_x(\bar{u}))) = f(k_x(\bar{u}), \bar{u}) = 0$ and therefore $\bar{x} := k_x(\bar{u})$ is an equilibrium of the closed-loop system. Conversely, let $\bar{x}$ be an equilibrium; the corresponding output value satisfies $\bar{y} = h(\bar{x})$. As in closed-loop $u = y$, we have $\bar{x} = k_x(\bar{y})$. Thus $\bar{y} = h(k_x(\bar{y}))$, as desired. We verify next that the $\mathcal{L}^\infty$ induced gain of the linearized system (19), satisfies:
$$\gamma_\infty = k'_y(\bar{u}).$$

In fact, we have:
$$k'_y(u) = \frac{\partial h}{\partial x}(k_x(u))\, k'_x(u). \qquad (25)$$

Recalling that, by definition, $f(k_x(u), u) = 0$, we can compute the derivative of $k_x$ by differentiating:
$$\frac{\partial f}{\partial x}(k_x(u), u)\, k'_x(u) + \frac{\partial f}{\partial u}(k_x(u), u) = 0 \qquad (26)$$

(at those points where the derivative exist). Evaluating the above expression at $u = \bar{u}$ yields $k'_x(\bar{u}) = -A^{-1}B$, where $A$ and $B$ are defined as $A = \left.\frac{\partial f}{\partial x}\right|_{x=k_x(\bar{u}), u=\bar{u}}$ and $B = \left.\frac{\partial f}{\partial u}\right|_{x=k_x(\bar{u}), u=\bar{u}}$ and $A^{-1}$ exists by non-degeneracy of the I/S characteristic. The claim follows by Remark 3.2. Next we investigate stability of the linearized system by looking at the slopes of the I/O characteristics at the equilibria.

We now make some remarks concerning the closed-loop linearized system, viz. the system arising by linearizing the nonlinear system (24) together with the unitary feedback interconnection $u = y$ (this is actually the same system resulting by first linearizing (24) and then applying unitary feedback). This system satisfies the equations $\dot{z} = (A + BC)z$ with the notations introduced so far. Moreover, it is monotone by virtue of Lemma 3.3. Therefore, it admits a real dominant eigenvalue $\bar{\lambda}$, viz. an eigenvalue so that $\bar{\lambda} \geq \text{Re}[\lambda_i]$ for all $i = 1, 2, \ldots n$ and $\bar{v} \in K$ for a corresponding eigenvector: $(A+BC)\bar{v} = \bar{\lambda}\bar{v}$. If we multiply both sides of the equality times $CA^{-1}$ we obtain:
$$\bar{\lambda}(CA^{-1}\bar{v}) = (C\bar{v})[1 + CA^{-1}B] = (C\bar{v})[1 - k'_y(\bar{u})] \qquad (27)$$

Moreover, by the asymptotic stability assumption on $A$:
$$CA^{-1}\bar{v} = -\int_0^{+\infty} C\underbrace{e^{At}\bar{v}}_{\in K}\, dt < 0, \qquad (28)$$

where the integral in (28) converges as $A$ is Hurwitz by virtue of Lemma 3.4.

Thus $\bar{\lambda} < 0$ if and only if $k'_y(\bar{u}) < 1$. On the other hand, $k'_y(\bar{u}) > 1$ iff $\bar{\lambda} > 0$. In particular, equilibria with $k'_y(\bar{u}) < 1$ are locally asymptotically stable and equilibria with $k'_y(\bar{u}) > 1$ have a nontrivial unstable manifolds. It is a routine exercise in measure theory to prove that the set of initial conditions which result in a converging trajectory is indeed a Borel set. Therefore, by Hirsch's Theorem on generic convergence of strongly monotone flows (see [9], Section 7) for almost all initial conditions, solutions will converge to the set of equilibria. Moreover, by Remark (4.3) the stable manifolds of (exponentially) unstable equilibria have zero-measure. Therefore, for almost all initial conditions solutions converge to points where $k'_y(\bar{u}) < 1$. This completes the proof of our result. ∎

**Remark 4.1** It is worth pointing out that, whenever the equilibrium in Theorem 7 is unique, convergence to the equilibrium is global under mild approximability assumptions which are always satisfied, for instance, if the state-space is convex. This is proved in Theorem 3.1 of [13]. □

**Remark 4.2** An alternative proof, based on frequency domain considerations, of the connection between stability of the closed-loop equilibrium and the I/O characteristic is provided next.

Consider the transfer function
$$w(s) = \int_0^\infty h(t) e^{-st} dt$$
of a strictly proper system, and let
$$w_{cl}(s) = \frac{w(s)}{1 - w(s)}$$
be the transfer function of the associated unity-feedback closed-loop system. Suppose:

1. $h$ is integrable (so, $w$ has no real nonnegative poles);
2. $h(t) \geq 0$ for all $t \geq 0$ (and is not identically zero);
3. $w(0) \neq 0$ (transversality condition).

Then:

(a) there exists a positive real pole of $w_{cl}$ if and only if $w(0) > 1$;

(b) every real pole of $w_{cl}$ is negative if and only if $w(0) < 1$.

*Proof:*

By the first assumption, $w(\lambda)$ is a continuous (real-valued) function for $\lambda \geq 0$.

Furthermore, $h(t) \geq 0$ for all $t \geq 0$ implies that $w'(\lambda) = -\int_0^\infty h(t) t e^{-\lambda t} dt < 0$ for all $\lambda$, so $w$ is a strictly decreasing function of $\lambda$.

Positive real poles of $w_{cl}$ are exactly those $\lambda > 0$ such that $w(\lambda) = 1$.

If $w(\lambda) = 1$ for some $\lambda > 0$ then the strict decrease of $w$ implies that $w(0) > 1$. Conversely, suppose that $w(0) > 1$. By strict properness, $w(\lambda) \to 0$ as $\lambda \to +\infty$. Thus there is some $\lambda > 0$ such that $w(\lambda) = 1$. This proves (a).

The first conclusion may be restated as: "every pole of $w_{cl}$ is $\leq 0$ if and only if $w(0) \leq 1$" so, since we know in addition that $w(0) \neq 1$, this is the same as requiring that every real pole is (strictly) negative. Thus (b) holds too.

□

**Remark 4.3** Stable manifolds of (exponentially) unstable equilibria have zero-measure. In the non-necessarily hyperbolic case, this fact is an easy consequence of Theorem 2.1 in [4] (modified as discussed in the remarks following Theorem 2.1, including the choice of suitable norms and the multiplication by a "bump" function, after a linear change of coordinates, and specialized to $r = 1$, and applied to time-1 maps). □

**Remark 4.4** A precise characterization of the basin of attraction of each asymptotically stable equilibrium is of course not possible in general; on the other hand, it is a straightforward consequence of monotonicity of the $I/S$ characteristic that equilibria are ordered, $e_1 \prec e_2 \prec e_3$. It therefore makes sense to speak about intervals $[e_1, e_2] := \{x \in X : e_1 \preceq x \preceq e_2\}$. Again, it is a straightforward consequence of monotonicity that intervals $[e_1, e_2]$ with $e_1, e_2$ equilibria are positively invariant. This allows to give estimates of the basin of attraction of each equilibrium. In the case of 3 equilibria for instance, with $e_1 \prec e_2 \prec e_3$ and $e_1, e_3$ asymptotically stable, $e_2$ unstable, we can conclude that $\{x : x \ll e_2\} \subset \mathcal{A}_1$ and $\{x : x \gg e_2\} \subset \mathcal{A}_3$. Similar considerations, based on empirical evidence, are made for instance in [3]. It is therein pointed-out how the unstable equilibrium plays the role of a threshold. □

## 5 Applications

A typical situation for the application of Theorem 7 is when a monotone system with a well-defined I/O characteristic of sigmoidal shape is closed under unitary feedback. If the sigmoidal function is sufficiently steep, this configuration is known to yield 3 equilibria, 2 stable and 1 unstable. In biological examples this might arise when a feedback loop comprising any number of positive interactions and an even number of inhibitions is present (no inhibition at all is also a situation which might lead to the same type of behavior). This is a well-known principle in biology. One of its simplest manifestations is the so called "competitive exclusion" principle, in which two competing species are coexisting and the possible equilibria are those where either one of the species is strongly inhibited. Just as an example consider the system described in [8, 11], describing the synthesis of lactose operon in Escherichia coli.. The systems equations are as follows:

$$\begin{aligned} \dot{x}_1 &= \frac{\alpha_1}{1+x_2^\beta} - x_1 \\ \dot{x}_2 &= \frac{\alpha_2}{1+x_1^\gamma} - x_2 \end{aligned} \qquad (29)$$

This can be seen as the unitary feedback interconnection of:

$$\begin{aligned} \dot{x}_1 &= \frac{\alpha_1}{1+u^\beta} - x_1 \\ \dot{x}_2 &= \frac{\alpha_2}{1+x_1^\gamma} - x_2 \end{aligned} \qquad (30)$$

$$y = x_2.$$

Equation (30) is a monotone dynamical system with respect to the order induced by the positivity cone $K := \mathbb{R}_{\leq 0} \times \mathbb{R}_{\geq 0}$. It is straightforward by a cascade argument to see that the system is endowed with the following static I/S characteristic:

$$k_x(u) = \begin{bmatrix} \frac{\alpha_1}{1+u^\beta} \\ \frac{\alpha_2(1+u^\beta)^\gamma}{(1+u^\beta)^\gamma + \alpha_1^\gamma} \end{bmatrix}.$$

In Fig. 1 we plotted the I/O static characteristic for $\alpha_1 = 1.3$, $\alpha_2 = 1$, $\beta = 3$ and $\gamma = 10$. As confirmed by simulations (see Phase plane), for almost all initial conditions trajectories converge to the equilibria where the derivative condition is satisfied. Of course, the interest of our results is in the high-dimensional case in which phase-plane techniques cannot provide the result.

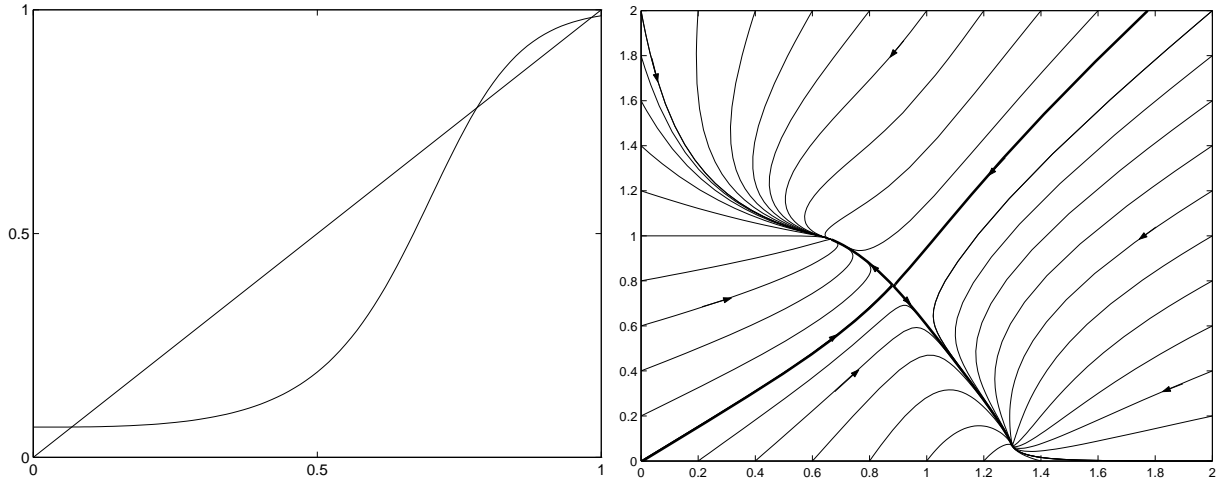

Figure 1: (a) I/O static characteristic; (b) Phase plane diagram

## 6  External Stimuli, Thresholds and Hysteresis

Throughout this section we investigate the behavior of positive feedback interconnections of monotone systems which are in turn excited by some exogenous input. In particular we consider interconnections of the following type:

$$\begin{aligned} \dot{x} &= f(x, u, v) \\ y &= h_y(x) \\ w &= h_w(x) \end{aligned} \qquad (31)$$

along with the unitary feedback interconnection $u = y$. The block diagram of such systems is shown in Fig. 2.

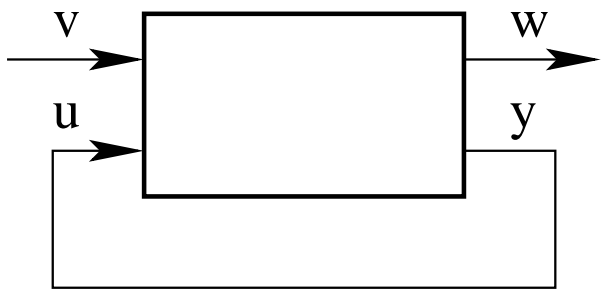

Figure 2: Block diagram of unitary feedback system with external inputs

We assume $f : X \times U \times V \to \mathbb{R}^n$ to be a locally Lipschitz function and that the system (31) is a monotone control system with input $[u, v]$ and output $[y, w]$ with respect to some ordering $\succeq_x$ of the state-space $X$ and cross-product orders as far as inputs $[u, v]$ and outputs $[y, w]$ are concerned, (viz. $[u_1, v_1] \succeq_I [u_2, v_2]$ iff $u_1 \succeq_u u_2$ and $v_1 \succeq_v v_2$, $[y_1, w_1] \succeq_O [y_2, w_2]$ iff $y_1 \succeq_y y_2$ and $w_1 \succeq_w w_2$).
For each fixed value of the input $v$, systems as in (31) can be studied according to the techniques described previously. A special instance of systems of this kind is given by SISO systems of the following form:

$$\begin{aligned} \dot{x} &= f(x, d) \\ d &= g(v, y) \\ y &= h(x) \end{aligned} \qquad (32)$$

where $g : V \times U \to \mathbb{R}$ is a monotone and locally Lipschitz function, (for instance $u, v \in \mathbb{R}_{\geq 0}$ and $g(v,y) = vy$ or $g(v,y) = v + y$ ). This structure ( see Fig. 3) is of interest because it arises commonly in biological applications and is particularly suited for a graphical analysis.

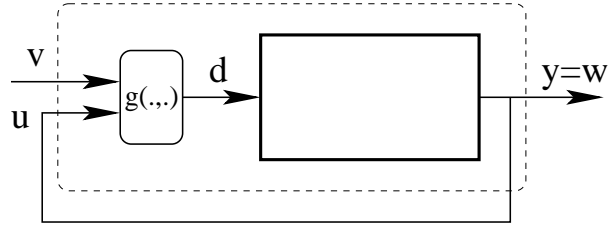

Figure 3: A special feedback configuration of SISO systems

Next, we discuss the behavior of such interconnections in the presence of external stimuli. In particular, in the case of multistable systems we prove the existence of *threshold* values of inputs which trigger the commutation between different equilibria.

The above considerations suggest the possibility of studying interconnections as in (31) by taking into account a parameterized family of I/O static characteristics in the $(u, y)$ plane, where the parameter is the exogenous input $v$. This type of analysis is very general and bifurcations can be traced by looking at the intersections of the parameterized I/O characteristic with the diagonal $u = y$. For the special structure (3) instead, the study can be carried out in the $(d, y)$-plane allowing some intuitive simplifications. A single I/O characteristic is needed in fact, from $d$ to $y$, and equilibria correspond to intersections with the "parameterized" family of functions $d = g(v, y)$, which also takes values in the $(d, y)$ plane. Although the analysis which follows is essentially a consequence of Theorem 7, it is still worth pursuing, because it provides a solid theoretical justification to phenomena which are well described and understood in many biological applications. Consider again the system (30), subject to the feedback interconnection $u = v \cdot y$. This results in the following set of equations:

$$
\begin{aligned}
\dot{x}_1 &= \frac{\alpha_1}{1+(v \cdot x_2)^\beta} - x_1 \\
\dot{x}_2 &= \frac{\alpha_2}{1+x_1^\gamma} - x_2
\end{aligned}
\tag{33}
$$

$$y = x_2.$$

We may therefore analyze the system by looking at the I/O static characteristic from $u$ to $y$, together with the $v$-parameterized family of lines $y = u/v$. Fig. 4 illustrates a typical situation, corresponding here to the parameters value in the following table:

| $\gamma$ | 6 | $\beta$ | 3 |
|---|---|---|---|
| $\alpha_1$ | 1.3 | $\alpha_2$ | 1.3 |

Notice that for $v = 1$ bistability is obtained; in particular two equilibria are asymptotically stable and one is an unstable saddle whose stable manifold behaves as a separatrix for the basins of attractions of the stable equilibria. Bifurcations occur at two different values of $v$, approximately $v_1 \approx 0.8$ and $v_2 \approx 1.35$. This values correspond to the slopes of the tangent lines to the I/O characteristic. For all $v > v_2$ in fact there only exists one equilibrium, usually referred to as the *activated* equilibrium. For $v < v_1$ again only one equilibrium occurs but corresponding to a *non-activated* state. These values play therefore the role of input *thresholds* that may trigger transition from the non-activated state to an activated one and vice versa. After a signal of amplitude bigger than $v_2$ is applied for a sufficiently long time, the state will be in proximity of the activated equilibria. Then, this level of output will be maintained even after $v(t)$ drops below $v_2$, provided that $v_1 < v(t)$. Further decrease of the $v(t)$ below $v_1$, for a sufficiently long time, will instead trigger transition to a deactivated state, which is afterward maintained also for higher values of $v(t)$, provided that $v(t) < v_2$. This kind of behavior, known as hysteresis, has been observed in many biological systems ( see for instance [7]). These techniques hence provide a way of computing threshold values even without knowing explicitly the systems equations or parameters, if the function $g(\cdot,\cdot)$ in (32) is known together with the I/O static characteristic. In other words many interesting deduction on the behavior of the system are possible only by checking (experimentally or in simulations) the static properties of the system.

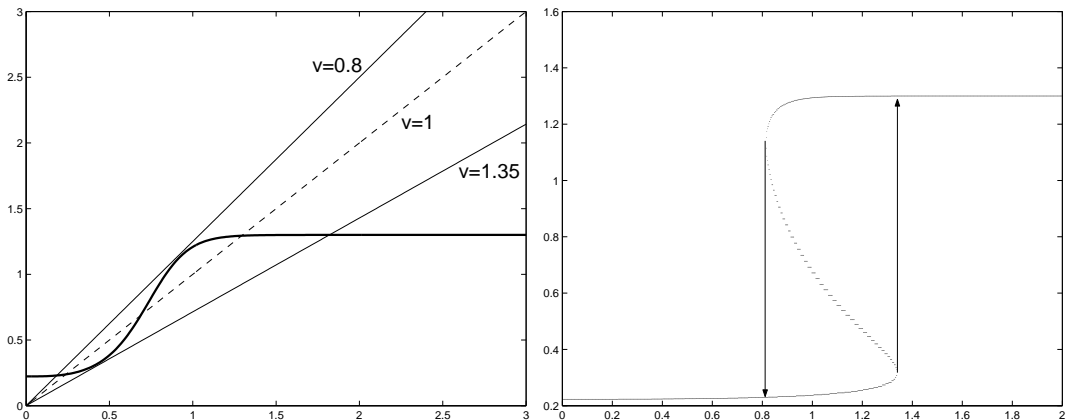

Figure 4: Thresholds and hysteresis

## 7 Why is Monotonicity Needed ?

Local analysis techniques based on the study of intersections of static characteristics of interconnected systems or, in the two-dimensional case of nullclines, are very common in mathematical biology. Our discussion shows that for the class of monotone systems, under relatively mild assumptions, almost global convergence results can be obtained and the investigation of the stability property of equilibria can be carried out just by graphical inspection at the intersection points of the I/O characteristics of systems in feedback. In this section we show by means of an example how monotonicity is a crucial assumption in this respect. The following planar system:

$$\begin{array}{rcl} \dot{x}_1 &=& x_1(-x_1 + x_2) \\ \dot{x}_2 &=& 3x_2(-x_1 + u) \\ y &=& c + b\frac{x_2^4}{k+x_2^4} \end{array} \qquad (34)$$

evolving in $\mathbb{R}^2_{\geq 0}$, it is not monotone. However, it has a well defined (monotonically increasing) I/O static characteristic, provided that $c, b, k \in \mathbb{R}_{>0}$. Moreover, for certain parameters values, the I/O characteristics has 3 (non-degenerate) fixed points. The closed-loop system resulting from the interconnection $u = y$, however, need not be globally converging at the set of equilibria. The simulations in Fig. 5 refer to the following values: $c = 1.1, b = 361/140, k = 405/14$. Notice that the 3 equilibria correspond to 2 unstable foci and one saddle point.

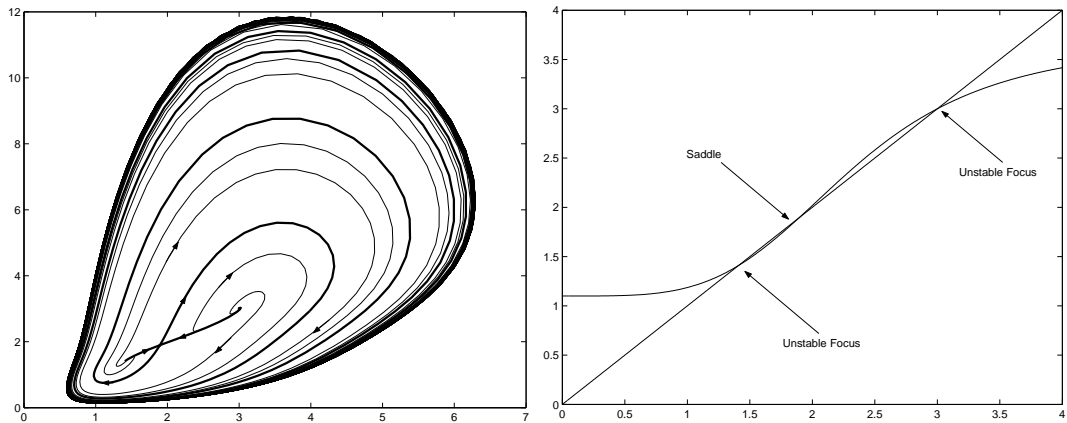

Figure 5: A stable limit cycle arising in a non-monotone feedback loop and the I/O static characteristic

## A  Graphical characterizations of Transparency and Excitability

**Proof of Theorem 2** Consider a monotone system which admits an incidence graph. Let $\xi$ be an arbitrary initial condition and $u(t)$, $v(t)$ arbitrary input signals satisfying $v(t) \gg u(t)$ for almost all $t \geq 0$. Clearly $x(t, \xi, v) \succeq x(t, \xi, u)$. Hence, we only need to show that inequality is actually strict for state components, viz. $x(t, \xi, v) \gg x(t, \xi, u)$. The result is proven by induction, by exploiting repeatedly Lemma 2.4. To this end it is useful to introduce decompositions of the system in sublayers, based on the following notion of distance between vertices of a graph:

$$d(v \to w) = \min\{L(\mathcal{P}) : \mathcal{P}_0 = v \text{ and } \mathcal{P}_L = w\}, \tag{35}$$

i.e. $d(v, w)$ denotes the shortest length among paths which link $v$ to $w$. For each vertex $x_i$ of the incidence graph corresponding to one of the state variable we define the integer:

$$D(x_i) = \min_j \; d(u_j \to x_i), \tag{36}$$

in words, $D(x_i)$ corresponds to the minimum distance from some input vertex $u_j$ to the state vertex $x_i$. By assumption $D(x_i)$ is well-defined for all $i \in \{1, \dots n\}$. We say that $x_i$ belongs to the $i$-th sublayer, if $D(x_i) = 1$. Consider any state component $x_i$ so that $D(x_i) = 1$ (by assumption such a component always exists). Clearly $\dot{x}_i = f_i(x, u)$ can be seen as a scalar (monotone) system, forced by the inputs $u$ and $x_{j \neq i}$. As $D(x_i) = u$, there exists $j$ such that $f$ is strictly monotone as a function of $u_j$ and therefore, $v_j(t) > u_j(t)$ for almost all $t \geq 0$ implies $x_i(t, \xi, v) > x_i(t, \xi, u)$ for all $t > 0$. By induction, any component belonging to the $i$-th sublayer is at least reachable in one step by some component belonging to the $i-1$-th sublayer, and therefore a similar argument applies, yielding $x_i(t, \xi, v) > x_i(t, \xi, u)$ for all $t > 0$.

This completes the proof for the case of weak excitability. Next we provide a detailed argument for the case of excitability. Assume that $v(t) \succ u(t)$ for almost all $t \geq 0$. Then, there exists an integer $j^\star$ so that $\{t : v_j(t) > u_j(t)\} \cap [0, \varepsilon)$ has non-zero measure for all $\varepsilon > 0$. We prove the result by induction by considering a sublayer decomposition taken by looking at graph distances with respect to the input vertex $j^\star$, viz. $D(x_i) := d(u_{j^\star} \to x_i)$. By assumption $D(x_i)$ is well-defined for all $i \in 1, 2, \dots n$. Consider any state component $x_i$ so that $D(x_i) = 1$ (such a component always exists). Again $\dot{x}_i = f_i(x, u)$ can be seen as a scalar (monotone) system, forced by the inputs $u$ and $x_{j \neq i}$. In particular $f_i(x, u)$ is strictly monotone with respect to $j^\star$. We want to show that $x_i(t, \xi, v) > x_i(t, \xi, u)$ for all $t > 0$. Assume by contradiction $x_i(t, \xi, v) = x_i(t, \xi, u) := x(t)$ for all $t \in [0, \varepsilon)$, for some $\varepsilon > 0$. Taking derivatives we have $f_i(x(t), v(t)) = f_i(x(t), u(t))$ for all $t \in [0, \varepsilon)$, but this contradicts strict monotonicity with respect to $u_{j^\star}$. By induction, any component belonging to the $i$-th sublayer is at least reachable in one step by some component belonging to the $i-1$-th sublayer, and therefore a similar argument applies, yielding $x_i(t, \xi, v) > x_i(t, \xi, u)$ for all $t > 0$.

**Proof of Theorem 3** Consider an arbitrary pair of ordered initial conditions $\xi_1 \succ \xi_2$. By monotonicity and unicity of solutions we have $x(t, \xi_1, u) \succ x(t, \xi_2, u)$ for all $t \geq 0$. As $x$ is finite-dimensional some $j^\star$ exists so that $\{t : x_{j^\star}(t, \xi_1, u) > x_{j^\star}(t, \xi_2, u)\} \cap [0, \varepsilon)$ has non-zero measure for all $\varepsilon > 0$. We claim that $x_i(t, \xi_1, u) > x_i(t, \xi_2, u)$ for all vertex $x_i$ reachable from the $x_{j^\star}$ and denote with $\mathcal{R}_{j^\star}$ the set of such $x_i$s. The claim can be shown inductively by arguments analogous to the one employed in the proof of Theorem 2.

By the graph reachability condition, for all (some) output vertices $y_j$ there exists at least $x_i \in \mathcal{R}_{j^\star}$ so that $x_i \to y_j$ is an edge of the incidence graph. Thus, $h_j(x(t, \xi_1, u)) > h_j(x(t, \xi_2, u))$ for all $t > 0$ for all such $j$s. This concludes the proof of the Theorem.